\documentclass[a4paper,english]{amsart}
\usepackage[T1]{fontenc}
\usepackage[latin9]{inputenc}
\usepackage[active]{srcltx}
\usepackage{amsthm}
\usepackage{amstext}
\usepackage{amssymb}

\makeatletter

\numberwithin{equation}{section}
\numberwithin{figure}{section}
\theoremstyle{plain}
\newtheorem{thm}{\protect\theoremname}
  \theoremstyle{plain}
  \newtheorem{prop}[thm]{\protect\propositionname}
  \theoremstyle{remark}
  \newtheorem*{rem*}{\protect\remarkname}

\makeatother

\usepackage{babel}
  \providecommand{\propositionname}{Proposition}
  \providecommand{\remarkname}{Remark}
\providecommand{\theoremname}{Theorem}

\begin{document}

\title{Poisson-Charlier and poly-Cauchy mixed type polynomials}

\author{By\\
 Dae San Kim And Taekyun Kim}
\begin{abstract}
In this paper, we consider Poisson-Charlier and poly-Cauchy mixed
type polynomials and give various identities of those polynomials
which are derived from umbral calculus.
\end{abstract}
\maketitle

\section{Introduction and Preliminaries}

$\,$

\global\long\def\lif{\textnormal{Lif}}

For $r\in\mathbb{Z}_{\ge0}$, the Cauchy numbers of the first kind
with order $r$ are defined by the generating function to be

\begin{center}
\begin{equation}
{\displaystyle \left(\frac{t}{\log\left(1+t\right)}\right)^{r}=\sum_{n=0}^{\infty}\mathbb{C}_{n}^{\left(r\right)}\frac{t^{n}}{n!},}\:\:(\textrm{see}\,[3,10,11,12]).\label{eq:1}
\end{equation}

\par\end{center}

In particular, when $r=1$, $\mathbb{C}_{n}^{\left(1\right)}=C_{n}$
are called Cauchy numbers of the first kind.

The Cauchy numbers of the second kind with order $r$ are defined
by
\begin{equation}
\left(\frac{t}{\left(1+t\right)\log\left(1+t\right)}\right)^{r}=\sum_{n=0}^{\infty}\hat{\mathbb{C}}_{n}^{\left(r\right)}\frac{t^{n}}{n!},\,\,(\textrm{see}\,[3,10,11,12]).\label{eq:2}
\end{equation}

When $r=1$, $\mathbb{\hat{C}}_{n}^{\left(1\right)}=\hat{C_{n}}$
are called the Cauchy numbers of the second kind.

As is well known, the generating function for the Poisson-Charlier
polynomials is given by
\begin{equation}
e^{-t}\left(1+\frac{t}{a}\right)^{x}=\sum_{n=0}^{\infty}C_{n}\left(x:a\right)\frac{t^{n}}{n!},\quad\left(a\ne0\right),\:(\textrm{see [14, 15])}.\label{eq:3}
\end{equation}

Recently, Komatsu has considered the poly-Cauchy polynomials of the
first kind as follows :
\begin{equation}
\frac{1}{\left(1+t\right)^{x}}\lif_{k}\left(\log\left(1+t\right)\right)=\sum_{n=0}^{\infty}C_{n}^{\left(k\right)}\left(x\right)\frac{t^{n}}{n!},\label{eq:4}
\end{equation}
where
\begin{equation}
\lif_{k}\left(t\right)=\sum_{n=0}^{\infty}\frac{t^{n}}{n!\left(n+1\right)^{k}},\:(\textrm{see [10, 11]}).\label{eq:5}
\end{equation}

He also introduced the poly-Cauchy polynomials of the second kind
by
\begin{equation}
\left(1+t\right)^{x}\lif_{k}\left(-\log\left(1+t\right)\right)=\sum_{n=0}^{\infty}\hat{C}_{n}^{\left(k\right)}\left(x\right)\frac{t^{n}}{n!},\:(\textrm{see [10, 11]}).\label{eq:6}
\end{equation}

In this paper, we consider Poisson-Charlier and poly-Cauchy of the
first kind mixed type polynomials as follows :
\begin{equation}
e^{-t}\lif_{k}\left(\log\left(1+\frac{t}{a}\right)\right)\left(1+\frac{t}{a}\right)^{-x}=\sum_{n=0}^{\infty}PC_{n}^{\left(k\right)}\left(x:a\right)\frac{t^{n}}{n!},\;\left(a\ne0\right).\label{eq:7}
\end{equation}

The Poisson-Charlier and poly-Cauchy of the second kind mixed type
polynomials are defined by the generating function to be
\begin{equation}
e^{-t}\lif_{k}\left(-\log\left(1+\frac{t}{a}\right)\right)\left(1+\frac{t}{a}\right)^{x}=\sum_{n=0}^{\infty}P\hat{C}_{n}^{\left(k\right)}\left(x:a\right)\frac{t^{n}}{n!},\;\left(a\ne0\right).\label{eq:8}
\end{equation}

It is known that the Frobenius-Euler polynomials of order $r$ are
given by
\begin{equation}
\left(\frac{1-\lambda}{e^{t}-\lambda}\right)^{r}e^{xt}=\sum_{n=0}^{\infty}H_{n}^{\left(r\right)}\left(x|\lambda\right)\frac{t^{n}}{n!},\:\left(\textrm{see [1, 4, 7, 9]}\right),\label{eq:9}
\end{equation}
where $r\in\mathbb{Z}_{\ge0}$, and $\lambda\in\mathbb{C}$ with $\lambda\ne1$.

The Bernoulli polynomials of order $r$ are also defined by the generating
function to be
\begin{equation}
\left(\frac{t}{e^{t}-1}\right)^{r}e^{xt}=\sum_{n=0}^{\infty}B_{n}^{\left(r\right)}\left(x\right)\frac{t^{n}}{n!},\:\left(\textrm{see [2, 5, 9, 10, 13]}\right).\label{eq:10}
\end{equation}
The Stirling number of the first kind is given by
\begin{equation}
\left(x\right)_{n}=\left(x\right)\left(x-1\right)\cdots\left(x-n+1\right)=\sum_{l=0}^{n}S_{1}\left(n,\, l\right)x^{l},\:\left(\textrm{see [14, 15]}\right),\label{eq:11}
\end{equation}
and by (\ref{eq:11}), we get
\begin{equation}
\left(\log\left(1+t\right)\right)^{m}=m!\sum_{l=m}^{\infty}S_{1}\left(l,\, m\right)\frac{t^{l}}{l!},\;\left(\textrm{see [8, 9, 14, 15]}\right).\label{eq:12}
\end{equation}

From (\ref{eq:11}), we note that
\begin{equation}
x^{\left(n\right)}=\left(-1\right)^{n}\left(-x\right)_{n}=\sum_{l=0}^{n}\left(-1\right)^{n-l}S_{1}\left(n,\, l\right)x^{l},\label{eq:13}
\end{equation}
where $x^{\left(n\right)}=x\left(x+1\right)\cdots\left(x+n-1\right),$
(see {[}1-15{]}).

Let $\mathbb{C}$ be the complex number field and let $\mathcal{F}$
be the set of all formal power series in the variable $t$ :
\begin{equation}
\mathcal{F}=\left\{ \left.f\left(t\right)=\sum_{k=0}^{\infty}\frac{a_{k}}{k!}t^{k}\right|a_{k}\in\mathbb{C}\right\} .\label{eq:14}
\end{equation}

Let $\mathbb{P}=\mathbb{C}\left[x\right]$ and let $\mathbb{P}^{*}$
be the vector space of all linear functionals on $\mathbb{P}$.

$\left\langle L|p\left(x\right)\right\rangle $ is the action of the
linear functional $L$ on the polynomial $p\left(x\right),$ and we
recall that the vector space operations on $\mathbb{P}^{*}$ are defined
by $\left\langle L+M|p\left(x\right)\right\rangle =\left\langle L|p\left(x\right)\right\rangle +\left\langle M|p\left(x\right)\right\rangle ,$
$\left\langle cL|p\left(x\right)\right\rangle =c\left\langle L|p\left(x\right)\right\rangle $,
where $c$ is complex constant in $\mathbb{C}$. For $f\left(t\right)\in\mathcal{F}$,
let us define the linear functional on $\mathbb{P}$ by setting
\begin{equation}
\left\langle f\left(t\right)|x^{n}\right\rangle =a_{n},\:\left(n\ge0\right).\label{eq:15}
\end{equation}
Thus, by (\ref{eq:14}) and (\ref{eq:15}), we get
\begin{equation}
\left\langle t^{k}|x^{n}\right\rangle =n!\:\delta_{n,k},\quad\left(n,\, k\ge0\right),\:\left(\textrm{see [4, 8, 14]}\right),\label{eq:16}
\end{equation}
where $\delta_{n,k}$ is the Kronecker's symbol.

Let ${\displaystyle f_{L}\left(t\right)=\sum_{k=0}^{\infty}\frac{\left\langle L|x^{k}\right\rangle }{k!}t^{k}}.$
Then, by (\ref{eq:15}), we see that $\left\langle \left.f_{L}\left(t\right)\right|x^{n}\right\rangle =\left\langle \left.L\right|x^{n}\right\rangle .$
The map $L\longmapsto f_{L}\left(t\right)$ is a vector space isomorphism
from $\mathbb{P}^{*}$ onto $\mathcal{F}$. Henceforth, $\mathcal{F}$
denotes both the algebra of formal power series in $t$ and the vector
space of all linear functionals on $\mathbb{P},$ and so an element
$f\left(t\right)$ of $\mathcal{F}$ will be thought of as both a
formal power series and a linear functional. We call $\mathcal{F}$
the umbral algebra and the umbral calculus is the study of umbral
algebra. The order $O\left(f\left(t\right)\right)$ of a power series
$f\left(t\right)\left(\ne0\right)$ is the smallest integer $k$ for
which the coefficient of $t^{k}$ does not vanish. If $O\left(f\left(t\right)\right)=1,$
then $f\left(t\right)$ is called a delta series; if $O\left(f\left(t\right)\right)=0,$
then $f\left(t\right)$ is called an invertible series. For $f\left(t\right),\, g\left(t\right)\in\mathcal{F}$
with $O\left(f\left(t\right)\right)=1$ and $O\left(g\left(t\right)\right)=0,$
there exists a unique sequence $s_{n}\left(x\right)$ such that $\left\langle \left.g\left(t\right)f\left(t\right)^{k}\right|s_{n}\left(x\right)\right\rangle =n!\delta_{n,k}$,
for $n,\, k\ge0.$ The sequence $s_{n}\left(x\right)$ is called the
sheffer sequence for $\left(g\left(t\right),\, f\left(t\right)\right)$
which is denoted by $s_{n}\left(x\right)\sim\left(g\left(t\right),\, f\left(t\right)\right)$
(see {[}8, 10, 14, 15{]}).

Let $f\left(t\right),\, g\left(t\right)\in\mathcal{F}$ and $p\left(x\right)\in\mathbb{P}$.
Then we see that
\begin{equation}
\left\langle \left.f\left(t\right)g\left(t\right)\right|p\left(x\right)\right\rangle =\left\langle \left.f\left(t\right)\right|g\left(t\right)p\left(x\right)\right\rangle =\left\langle \left.g\left(t\right)\right|f\left(t\right)p\left(x\right)\right\rangle ,\label{eq:17}
\end{equation}
and
\begin{equation}
f\left(t\right)=\sum_{k=0}^{\infty}\left\langle \left.f\left(t\right)\right|x^{k}\right\rangle \frac{t^{k}}{k!},\quad p\left(x\right)=\sum_{k=0}^{\infty}\left\langle \left.t^{k}\right|p\left(x\right)\right\rangle \frac{x^{k}}{k!}.\label{eq:18}
\end{equation}

By (\ref{eq:18}), we get
\begin{equation}
t^{k}p\left(x\right)=p^{\left(k\right)}\left(x\right)=\frac{d^{k}p\left(x\right)}{dx^{k}},\:\textrm{and }e^{yt}p\left(x\right)=p\left(x+y\right),\:\left(\textrm{see [14]}\right).\label{eq:19}
\end{equation}

For $s_{n}\left(x\right)\sim\left(g\left(t\right),\, f\left(t\right)\right),$
we have the generating function of $s_{n}\left(x\right)$ as follows
:
\begin{equation}
\frac{1}{g\left(\overline{f}\left(t\right)\right)}e^{x\overline{f}\left(t\right)}=\sum_{n=0}^{\infty}s_{n}\left(x\right)\frac{t^{n}}{n!},\textrm{ for all \ensuremath{x\in\mathbb{C},}}\label{eq:20}
\end{equation}
where $\overline{f}\left(t\right)$ is the compositional inverse of
$f\left(t\right)$ with $\overline{f}\left(f\left(t\right)\right)=t.$

Let $s_{n}\left(x\right)\sim\left(g\left(t\right),\, f\left(t\right)\right).$
Then we have the following equations (see {[}8, 14, 15{]}):

\begin{equation}
f\left(t\right)s_{n}\left(x\right)=ns_{n-1}\left(x\right),\;\left(n\ge0\right),\:\frac{d}{dx}s_{n}\left(x\right)=\sum_{l=0}^{n-1}\dbinom{n}{l}\left\langle \left.\overline{f}\left(t\right)\right|x^{n-l}\right\rangle s_{l}\left(x\right),\label{eq:21}
\end{equation}
\begin{equation}
s_{n}\left(x\right)=\sum_{j=0}^{n}\frac{1}{j!}\left\langle \left.g\left(\overline{f}\left(t\right)\right)^{-1}\overline{f}\left(t\right)^{j}\right|x^{n}\right\rangle x^{j},\:\left\langle \left.f\left(t\right)\right|xp\left(x\right)\right\rangle =\left\langle \left.\partial_{t}f\left(t\right)\right|p\left(x\right)\right\rangle ,\label{eq:22}
\end{equation}
and
\begin{equation}
s_{n}\left(x+y\right)=\sum_{j=0}^{n}\dbinom{n}{j}s_{j}\left(x\right)p_{n-j}\left(y\right),\:\textrm{where }p_{n}\left(x\right)=g\left(t\right)s_{n}\left(x\right).\label{eq:23}
\end{equation}

For $p_{n}\left(x\right)\sim\left(1,\, f\left(t\right)\right),$ $q_{n}\left(x\right)\sim\left(1,\, g\left(t\right)\right),$
it is well known that
\begin{equation}
q_{n}\left(x\right)=x\left(\frac{f\left(t\right)}{g\left(t\right)}\right)^{n}x^{-1}p_{n}\left(x\right),\;\left(n\ge1\right),\:\left(\textrm{see [14, 15]}\right).\label{eq:24}
\end{equation}

For $s_{n}\left(x\right)\sim\left(g\left(t\right),\, f\left(t\right)\right),$
$r_{n}\left(x\right)\sim\left(h\left(t\right),\, l\left(t\right)\right),$
let us assume that

\[
{\displaystyle s_{n}\left(x\right)=\sum_{m=0}^{\infty}C_{n,m}r_{n}\left(x\right),\mbox{\:}\left(n\ge0\right).}
\]

Then we have
\begin{equation}
C_{n,m}=\frac{1}{m!}\left\langle \left.\frac{h\left(\overline{f}\left(t\right)\right)}{g\left(\overline{f}\left(t\right)\right)}l\left(\overline{f}\left(t\right)\right)^{m}\right|x^{n}\right\rangle ,\quad\left(\textrm{see [8, 10, 14]}\right).\label{eq:25}
\end{equation}

In this paper, we investigate some identities of Poisson-Charlier
and poly-Cauchy mixed type polynomials arising from umbral calculus.
That is, we give various identities of the Poisson-Charlier and poly-Cauchy
polynomials of the first and second kind mixed type polynomials which
are derived from umbral calculus.

\section{Poisson-Charlier and poly-Cauchy mixed type polynomials}

$\,$

From (\ref{eq:6}), (\ref{eq:7}) and (\ref{eq:20}), we note that
\begin{equation}
PC_{n}^{\left(k\right)}\left(x:a\right)\sim\left(e^{a\left(e^{-t}-1\right)}\frac{1}{\lif_{k}\left(-t\right)},\, a\left(e^{-t}-1\right)\right),\label{eq:26}
\end{equation}
and

\begin{equation}
P\hat{C}_{n}^{\left(k\right)}\left(x:a\right)\sim\left(e^{a\left(e^{t}-1\right)}\frac{1}{\lif_{k}\left(-t\right)},\, a\left(e^{t}-1\right)\right).\label{eq:27}
\end{equation}

Now, we observe that

\begin{eqnarray}
 &  & PC_{n}^{\left(k\right)}\left(y:a\right)\label{eq:28}\\
 & = & \left\langle \left.\sum_{l=0}^{\infty}PC_{l}^{\left(k\right)}\left(y:a\right)\frac{t^{l}}{l!}\right|x^{n}\right\rangle \nonumber \\
 & = & \left\langle \left.e^{-t}\lif_{k}\left(\log\left(1+\frac{t}{a}\right)\right)\left(1+\frac{t}{a}\right)^{-y}\right|x^{n}\right\rangle \nonumber \\
 & = & \left\langle e^{-t}\left|\lif_{k}\left(\log\left(1+\frac{t}{a}\right)\right)\left(1+\frac{t}{a}\right)^{-y}x^{n}\right.\right\rangle \nonumber \\
 & = & \sum_{l=0}^{n}C_{l}^{\left(k\right)}\left(y\right)\frac{1}{a^{l}l!}\left(n\right)_{l}\left\langle \left.e^{-t}\right|x^{n-l}\right\rangle =\sum_{l=0}^{n}C_{l}^{\left(k\right)}\left(y\right)\dbinom{n}{l}\frac{\left(-1\right)^{n-l}}{a^{l}}.\nonumber
\end{eqnarray}

Therefore, by (\ref{eq:28}), we obtain the following theorem.
\begin{thm}
\label{thm:1}For $n\ge0,$ we have
\[
PC_{n}^{\left(k\right)}\left(x:a\right)=\sum_{l=0}^{n}C_{l}^{\left(k\right)}\left(x\right)\dbinom{n}{l}\frac{\left(-1\right)^{n-l}}{a^{l}},
\]
where $a\ne0$.
\end{thm}
$\,$

Alternatively,
\begin{eqnarray}
PC_{n}^{\left(k\right)}\left(y:a\right) & = & \left\langle \left.\lif_{k}\left(\log\left(1+\frac{t}{a}\right)\right)\right|e^{-t}\left(1+\frac{t}{a}\right)^{-y}x^{n}\right\rangle \label{eq:29}\\
 & = & \sum_{l=0}^{n}C_{l}\left(-y:a\right)\frac{1}{l!}\left(n\right)_{l}\left\langle \left.\lif_{k}\left(\log\left(1+\frac{t}{a}\right)\right)\right|x^{n-l}\right\rangle \nonumber \\
 & = & \sum_{l=0}^{n}C_{l}\left(-y:a\right)\dbinom{n}{l}C_{n-l}^{\left(k\right)}\frac{1}{a^{n-l}}\nonumber \\
 & = & \sum_{l=0}^{n}\frac{\begin{pmatrix}n\\
l
\end{pmatrix}C_{n-l}^{\left(k\right)}}{a^{n-l}}C_{l}\left(-y:a\right).\nonumber
\end{eqnarray}

Therefore, by (\ref{eq:29}), we obtain the following proposition.
\begin{prop}
\label{thm:2}For $n\ge0,$ $a\ne0$, we have
\[
PC_{n}^{\left(k\right)}\left(x:a\right)=\sum_{l=0}^{n}\frac{\begin{pmatrix}n\\
l
\end{pmatrix}C_{n-l}^{\left(k\right)}}{a^{n-l}}C_{l}\left(-x:a\right).
\]

\end{prop}
$\,$
\begin{rem*}
By the same method as (\ref{eq:28}) and (\ref{eq:29}), we get

\begin{equation}
P\hat{C}_{n}^{\left(k\right)}\left(x:a\right)={\displaystyle \sum_{l=0}^{n}\frac{\left(-1\right)^{n-l}\begin{pmatrix}n\\
l
\end{pmatrix}}{a^{l}}}\hat{C}_{l}^{\left(k\right)}\left(x\right),\label{eq:30}
\end{equation}

and

\begin{equation}
P\hat{C}_{n}^{\left(k\right)}\left(x:a\right)={\displaystyle \sum_{l=0}^{n}\frac{\begin{pmatrix}n\\
l
\end{pmatrix}\hat{C}_{l}^{\left(k\right)}}{a^{l}}C_{n-l}\left(x:a\right).}\label{eq:31}
\end{equation}

\end{rem*}
$\,$

It is not difficult to show that
\begin{equation}
\left(-\frac{1}{a}\right)^{n}x^{\left(n\right)}=a^{-n}\sum_{l=0}^{n}\left(-1\right)^{k}S_{1}\left(n,\, k\right)x^{k}\sim\left(1,\, a\left(e^{-t}-1\right)\right),\label{eq:32}
\end{equation}
and
\begin{equation}
a^{-n}\left(x\right)_{n}=a^{-n}\sum_{k=0}^{n}S_{1}\left(n,\, k\right)x^{k}\sim\left(1,\, a\left(e^{t}-1\right)\right).\label{eq:33}
\end{equation}

By (\ref{eq:26}), we get
\begin{equation}
e^{a\left(e^{-t}-1\right)}\frac{1}{\lif_{k}\left(-t\right)}PC_{n}^{\left(k\right)}\left(x:a\right)\sim\left(1,\, a\left(e^{-t}-1\right)\right).\label{eq:34}
\end{equation}
From (\ref{eq:32}), (\ref{eq:34}), we have
\begin{equation}
e^{a\left(e^{-t}-1\right)}\frac{1}{\lif_{k}\left(-t\right)}PC_{n}^{\left(k\right)}\left(x:a\right)=\left(-\frac{1}{a}\right)^{n}x^{\left(n\right)}.\label{eq:35}
\end{equation}
Thus, by (\ref{eq:35}) we get
\begin{eqnarray}
PC_{n}^{\left(k\right)}\left(x:a\right) & = & \lif_{k}\left(-t\right)e^{-a\left(e^{-t}-1\right)}\left(-\frac{1}{a}\right)^{n}x^{\left(n\right)}\label{eq:36}\\
 & = & \left(-\frac{1}{a}\right)^{n}\lif_{k}\left(-t\right)\sum_{l=0}^{n}\frac{a^{l}}{l!}\left(1-e^{-t}\right)^{l}x^{\left(n\right)}.\nonumber
\end{eqnarray}

By (\ref{eq:13}), we see that $x^{\left(n\right)}\sim\left(1,\,1-e^{-t}\right).$
From (\ref{eq:21}) and (\ref{eq:36}), we have
\[
\left(1-e^{-t}\right)^{l}x^{\left(n\right)}=\left(n\right)_{l}x^{\left(n-l\right)}.
\]
and
\begin{align}
 & PC_{n}^{\left(k\right)}\left(x:a\right)\label{eq:37}\\
= & \left(-\frac{1}{a}\right)^{n}\lif_{k}\left(-t\right)\sum_{l=0}^{n}\frac{a^{l}}{l!}\left(1-e^{-t}\right)^{l}x^{\left(n\right)}\nonumber \\
= & \left(-\frac{1}{a}\right)^{n}\lif_{k}\left(-t\right)\sum_{l=0}^{n}\begin{pmatrix}n\\
l
\end{pmatrix}a^{l}x^{\left(n-l\right)}\nonumber \\
= & \left(-\frac{1}{a}\right)^{n}\sum_{l=0}^{n}a^{l}\begin{pmatrix}n\\
l
\end{pmatrix}\sum_{m=0}^{n-l}\left(-1\right)^{n-l-m}S_{1}\left(n-l,\, m\right)\sum_{r=0}^{m}\frac{\left(-1\right)^{r}}{r!\left(r+1\right)^{k}}t^{r}x^{m}\nonumber \\
= & a^{-n}\sum_{l=0}^{n}\sum_{m=0}^{n-l}\sum_{r=0}^{m}\left(-1\right)^{l+m+r}\begin{pmatrix}n\\
l
\end{pmatrix}\begin{pmatrix}m\\
r
\end{pmatrix}\frac{a^{l}}{\left(r+1\right)^{k}}S_{1}\left(n-l,\, m\right)x^{m-r}\nonumber \\
= & a^{-n}\sum_{j=0}^{n}\left\{ \sum_{m=j}^{n}\sum_{l=0}^{n-m}\left(-1\right)^{l+j}\begin{pmatrix}n\\
l
\end{pmatrix}\begin{pmatrix}m\\
j
\end{pmatrix}\frac{a^{l}}{\left(m-j+1\right)^{k}}S_{1}\left(n-l,\, m\right)\right\} x^{j}.\nonumber
\end{align}

Therefore, by (\ref{eq:37}), we obtain the following theorem.
\begin{thm}
\label{thm:3}For $n\ge0$, $k\in\mathbb{Z}$ and $a\ne0$, we have
\begin{align*}
 & PC_{n}^{\left(k\right)}\left(x:a\right)\\
= & a^{-n}\sum_{j=0}^{n}\left\{ \sum_{m=j}^{n}\sum_{l=0}^{n-m}\left(-1\right)^{l+j}\begin{pmatrix}n\\
l
\end{pmatrix}\begin{pmatrix}m\\
j
\end{pmatrix}\frac{a^{l}}{\left(m-j+1\right)^{k}}S_{1}\left(n-l,\, m\right)\right\} x^{j}.
\end{align*}

\end{thm}
$\,$
\begin{rem*}
By (\ref{eq:27}) and (\ref{eq:33}), we get
\begin{eqnarray}
P\hat{C}_{n}^{\left(k\right)}\left(x:a\right) & = & \lif_{k}\left(-t\right)e^{-a\left(e^{t}-1\right)}a^{-n}\left(x\right)_{n}\label{eq:38}\\
 & = & a^{-n}\lif_{k}\left(-t\right)\sum_{l=0}^{n}\frac{\left(-a\right)^{l}}{l!}\left(e^{t}-1\right)^{l}\left(x\right)_{n}.\nonumber
\end{eqnarray}

\end{rem*}
$\,$

By the same method as (\ref{eq:37}), we get
\begin{align*}
 & P\hat{C}_{n}^{\left(k\right)}\left(x:a\right)\\
= & a^{-n}\sum_{j=0}^{n}\left\{ \sum_{m=j}^{n}\sum_{l=0}^{n-m}\left(-1\right)^{l+m+j}\begin{pmatrix}n\\
l
\end{pmatrix}\begin{pmatrix}m\\
j
\end{pmatrix}\frac{a^{l}}{\left(m-j+1\right)^{k}}S_{1}\left(n-l,\, m\right)\right\} x^{j}.
\end{align*}

From (\ref{eq:22}) and (\ref{eq:26}), we note that
\begin{align}
 & PC_{n}^{\left(k\right)}\left(x:a\right)\label{eq:40}\\
= & \sum_{l=0}^{n}\frac{1}{l!}\left\langle \left.e^{-t}\lif_{k}\left(\log\left(1+\frac{t}{a}\right)\right)\left(-\log\left(1+\frac{t}{a}\right)\right)^{l}\right|x^{n}\right\rangle x^{l}\nonumber \\
= & \sum_{l=0}^{n}\frac{\left(-1\right)^{l}}{l!}\sum_{r=0}^{n-l}\frac{l!}{\left(r+l\right)!a^{r+l}}S_{1}\left(r+l,\, l\right)\left\langle \left.e^{-t}\lif_{k}\left(\log\left(1+\frac{t}{a}\right)\right)\right|t^{r+n}x^{n}\right\rangle x^{l}\nonumber \\
= & \sum_{l=0}^{n}\left\{ \sum_{r=0}^{n-l}\frac{\left(-1\right)^{l}}{a^{r+l}}\begin{pmatrix}n\\
r+l
\end{pmatrix}S_{1}\left(r+l,\, l\right)PC_{n-r-l}^{\left(k\right)}\left(0:a\right)\right\} x^{l}\nonumber \\
= & \sum_{l=0}^{n}\left\{ \sum_{r=0}^{n-l}\frac{\left(-1\right)^{l}}{a^{n-r}}\begin{pmatrix}n\\
r
\end{pmatrix}S_{1}\left(n-r,\, l\right)PC_{r}^{\left(k\right)}\left(0:a\right)\right\} x^{l}.\nonumber
\end{align}

Therefore, by (\ref{eq:40}), we obtain the following theorem.
\begin{thm}
\label{thm:4} For $n\ge0$, $k\in\mathbb{Z}$ and $a\ne0$, we have
\[
PC_{n}^{\left(k\right)}\left(x:a\right)=\sum_{l=0}^{n}\left\{ \sum_{r=0}^{n-l}\frac{\left(-1\right)^{l}}{a^{n-r}}\begin{pmatrix}n\\
r
\end{pmatrix}S_{1}\left(n-r,\, l\right)PC_{r}^{\left(k\right)}\left(0:a\right)\right\} x^{l}.
\]
\end{thm}
\begin{rem*}
From (\ref{eq:22}) and (\ref{eq:27}), we can also derive the following
equation.

\begin{equation}
P\hat{C}_{n}^{\left(k\right)}\left(x:a\right)=\sum_{l=0}^{n}\left\{ \sum_{r=0}^{n-l}\frac{\begin{pmatrix}n\\
r
\end{pmatrix}}{a^{n-r}}S_{1}\left(n-r,\, l\right)P\hat{C_{r}}^{\left(k\right)}\left(0:a\right)\right\} x^{l}.\label{eq:41}
\end{equation}

\end{rem*}
$\,$

By (\ref{eq:26}), we easily see that
\begin{equation}
e^{a\left(e^{-t}-1\right)}\frac{1}{\lif_{k}\left(-t\right)}PC_{n}^{\left(k\right)}\left(x:a\right)\sim\left(1,\, a\left(e^{-t}-1\right)\right),\quad x^{n}\sim\left(1,\, t\right).\label{eq:42}
\end{equation}

Thus, by (\ref{eq:24}) and (\ref{eq:42}), for $n\ge1$ we get
\begin{align}
 & e^{a\left(e^{-t}-1\right)}\frac{1}{\lif_{k}\left(-t\right)}PC_{n}^{\left(k\right)}\left(x:a\right)\label{eq:43}\\
= & x\left(\frac{t}{a\left(e^{-t}-1\right)}\right)^{n}x^{-1}x^{n}=\left(-a^{-1}\right)^{n}x\left(\frac{-t}{e^{-t}-1}\right)^{n}x^{n-1}\nonumber \\
= & \left(-a^{-1}\right)^{n}x\sum_{r=0}^{\infty}B_{r}^{\left(n\right)}\frac{\left(-t\right)^{r}}{r!}x^{n-1}\nonumber \\
= & \left(-a^{-1}\right)^{n}x\sum_{r=0}^{n-1}B_{r}^{\left(n\right)}\left(n-1\right)_{r}\frac{\left(-1\right)^{r}}{r!}x^{n-r-1}\nonumber \\
= & \left(-a^{-1}\right)^{n}\sum_{r=0}^{n-1}B_{r}^{\left(n\right)}\left(-1\right)^{r}\begin{pmatrix}n-1\\
r
\end{pmatrix}x^{n-r},\nonumber
\end{align}
where $B_{r}^{\left(n\right)}=B_{r}^{\left(n\right)}\left(0\right)$
are called the Bernoulli numbers of order $n$.

From (\ref{eq:43}), we have
\begin{align}
 & PC_{n}^{\left(k\right)}\left(x:a\right)\label{eq:44}\\
= & \left(-a^{-1}\right)^{n}\sum_{r=0}^{n-1}\left(-1\right)^{r}\begin{pmatrix}n-1\\
r
\end{pmatrix}B_{r}^{\left(n\right)}\lif_{k}\left(-t\right)e^{-a\left(e^{-t}-1\right)}x^{n-r}\nonumber \\
= & \left(-a^{-1}\right)^{n}\sum_{r=0}^{n-1}\left(-1\right)^{r}\begin{pmatrix}n-1\\
r
\end{pmatrix}B_{r}^{\left(n\right)}\lif_{k}\left(-t\right)\sum_{l=0}^{n-r}\frac{\left(-a\right)^{l}}{l!}\left(e^{-t}-1\right)^{l}x^{n-r}\nonumber \\
= & \left(-a^{-1}\right)^{n}\sum_{r=0}^{n-1}\left(-1\right)^{r}\begin{pmatrix}n-1\\
r
\end{pmatrix}B_{r}^{\left(n\right)}\sum_{l=0}^{n-r}\left\{ \frac{\left(-a\right)^{l}}{l!}\sum_{j=0}^{n-r-l}\frac{l!}{\left(j+l\right)!}S_{2}\left(j+l,\, l\right)\right.\nonumber \\
 & \left.\qquad\qquad\quad\times\left(-1\right)^{j+l}\lif_{k}\left(-t\right)t^{j+l}x^{n-r}\right\} ,\nonumber
\end{align}
where $S_{2}\left(n,\, k\right)$ is the stirling number of the second
kind.

Now, we observe that
\begin{align}
 & \lif_{k}\left(-t\right)t^{j+l}x^{n-r}\label{eq:45}\\
= & \left(n-r\right)_{j+l}\lif_{k}\left(-t\right)x^{n-r-j-l}\nonumber \\
= & \left(n-r\right)_{j+l}\sum_{m=0}^{\infty}\frac{\left(-1\right)^{m}t^{m}}{m!\left(m+1\right)^{k}}x^{n-r-j-l}\nonumber \\
= & \left(n-r\right)_{j+l}\sum_{m=0}^{n-r-j-l}\frac{\left(-1\right)^{m}}{m!\left(m+1\right)^{k}}\left(n-r-j-l\right)_{m}x^{n-r-j-l-m}.\nonumber
\end{align}

Therefore, by (\ref{eq:44}) and (\ref{eq:45}), we obtain the following
theorem.
\begin{thm}
\label{thm:5} For $n\ge1$, $k\in\mathbb{Z}$ and $a\ne0$, we have

\begin{align*}
 & PC_{n}^{\left(k\right)}\left(x:a\right)\\
= & a^{-n}\sum_{m=0}^{n}\left\{ \sum_{r=0}^{n-m}\sum_{l=0}^{n-m-r}\sum_{j=0}^{n-m-r-l}\left(-1\right)^{l+m}\begin{pmatrix}n-1\\
r
\end{pmatrix}\begin{pmatrix}n-r\\
j+l
\end{pmatrix}\right.\\
 & \qquad\quad\:\left.\times\begin{pmatrix}n-r-j-l\\
m
\end{pmatrix}\frac{a^{l}S_{2}\left(j+l,\, l\right)}{\left(n-r-j-l-m+1\right)^{k}}B_{r}^{\left(n\right)}\right\} x^{m}.
\end{align*}

\end{thm}
$\,$
\begin{rem*}
We note that
\begin{equation}
e^{a\left(e^{t}-1\right)}\frac{1}{\lif_{k}\left(-t\right)}P\hat{C}_{n}^{\left(k\right)}\left(x:a\right)\sim\left(1,a\left(e^{t}-1\right)\right),\quad x^{n}\sim\left(1,\, t\right).\label{eq:46}
\end{equation}

Thus, for $n\ge1$ we have
\begin{equation}
e^{a\left(e^{t}-1\right)}\frac{1}{\lif_{k}\left(-t\right)}P\hat{C}_{n}^{\left(k\right)}\left(x:a\right)=a^{-n}\sum_{l=0}^{n-1}\dbinom{n-1}{l}B_{l}^{\left(n\right)}x^{n-l}.\label{eq:47}
\end{equation}

\end{rem*}
$\,$

From (\ref{eq:47}), for $n\ge1$ we can derive
\begin{align}
 & P\hat{C}_{n}^{\left(k\right)}\left(x:a\right)\label{eq:48}\\
= & \left(-a^{-1}\right)^{n}\sum_{m=0}^{n}\left[\sum_{r=0}^{n-m}\sum_{l=0}^{n-m-r}\sum_{j=0}^{n-m-r-l}\left\{ \left(-1\right)^{r+j+m}\frac{\begin{pmatrix}n-1\\
r
\end{pmatrix}\begin{pmatrix}n-r\\
j+l
\end{pmatrix}}{\left(n-r-j-l-m+1\right)^{k}}\right.\right.\nonumber \\
 & \qquad\qquad\qquad\quad\left.\left.\times\begin{pmatrix}n-r-j-l\\
m
\end{pmatrix}a^{l}S_{2}\left(j+l,\, l\right)B_{r}^{\left(n\right)}\right\} \right]x^{m}.\nonumber
\end{align}

By (\ref{eq:23}), (\ref{eq:26}) and (\ref{eq:27}), we get

\begin{eqnarray}
PC_{n}^{\left(k\right)}\left(x+y:a\right) & = & \sum_{j=0}^{n}\dbinom{n}{j}PC_{j}^{\left(k\right)}\left(x:a\right)\left(-a^{-1}\right)^{n-j}y^{\left(n-j\right)}\label{eq:49}\\
 & = & \sum_{j=0}^{n}\dbinom{n}{j}PC_{n-j}^{\left(k\right)}\left(x:a\right)\left(-a^{-1}\right)^{j}y^{\left(j\right)}\nonumber
\end{eqnarray}
and
\begin{eqnarray}
P\hat{C}_{n}^{\left(k\right)}\left(x+y:a\right) & = & \sum_{j=0}^{n}\dbinom{n}{j}P\hat{C}_{j}^{\left(k\right)}\left(x:a\right)a^{-\left(n-j\right)}\left(y\right)_{n-j}\label{eq:50}\\
 & = & \sum_{j=0}^{n}\dbinom{n}{j}P\hat{C}_{n-j}^{\left(k\right)}\left(x:a\right)a^{-j}\left(y\right)_{j}.\nonumber
\end{eqnarray}

From (\ref{eq:21}), (\ref{eq:26}) and (\ref{eq:27}), we have
\begin{equation}
PC_{n}^{\left(k\right)}\left(x-1:a\right)-PC_{n}^{\left(k\right)}\left(x:a\right)=a^{-1}nPC_{n-1}^{\left(k\right)}\left(x:a\right),\label{eq:51}
\end{equation}
and
\begin{equation}
P\hat{C}_{n}^{\left(k\right)}\left(x+1:a\right)-PC_{n}^{\left(k\right)}\left(x:a\right)=a^{-1}nP\hat{C}_{n-1}^{\left(k\right)}\left(x:a\right).\label{eq:52}
\end{equation}

For $s_{n}\left(x\right)\sim\left(g\left(t\right),\, f\left(t\right)\right)$,
we note that recurrence formula for $s_{n}\left(x\right)$ is given
by
\begin{equation}
s_{n+1}\left(x\right)=\left(x-\frac{g^{\prime}\left(t\right)}{g\left(t\right)}\right)\frac{1}{f^{\prime}\left(t\right)}s_{n}\left(x\right).\label{eq:53}
\end{equation}

Thus, by (\ref{eq:26}), (\ref{eq:27}) and (\ref{eq:53}), we get
\begin{align}
 & PC_{n+1}^{\left(k\right)}\left(x:a\right)\label{eq:54}\\
= & -\frac{1}{a}xPC_{n}^{\left(k\right)}\left(x+1:a\right)-PC_{n}^{\left(k\right)}\left(x:a\right)\nonumber \\
 & +a^{-\left(n+1\right)}\sum_{j=0}^{n}\left\{ \sum_{m=j}^{n}\sum_{l=0}^{n-m}\left(-1\right)^{l+j}\dbinom{n}{l}\dbinom{m}{j}\frac{a^{l}}{\left(m-j+2\right)^{k}}S_{1}\left(n-l,m\right)\right\} \left(x+1\right)^{j},\nonumber
\end{align}
and
\begin{align}
 & P\hat{C}_{n+1}^{\left(k\right)}\left(x:a\right)\label{eq:55}\\
= & \frac{1}{a}xP\hat{C}_{n}^{\left(k\right)}\left(x-1:a\right)-P\hat{C}_{n}^{\left(k\right)}\left(x:a\right)\nonumber \\
 & -a^{-\left(n+1\right)}\sum_{j=0}^{n}\left\{ \sum_{m=j}^{n}\sum_{l=0}^{n-m}\left(-1\right)^{l+m+j}\dbinom{n}{l}\dbinom{m}{j}\right.\nonumber \\
 & \qquad\qquad\quad\quad\left.\times\frac{a^{l}}{\left(m-j+2\right)^{k}}S_{1}\left(n-l,m\right)\right\} \left(x-1\right)^{j}.\nonumber
\end{align}

Note that
\begin{align}
 & PC_{n}^{\left(k\right)}\left(y:a\right)\label{eq:56}\\
= & \left\langle \left.\sum_{l=0}^{\infty}PC_{l}^{\left(k\right)}\left(y:a\right)\frac{t^{l}}{l!}\right|x^{n}\right\rangle =\left\langle \left.e^{-t}\lif_{k}\left(\log\left(1+\frac{t}{a}\right)\right)\left(1+\frac{t}{a}\right)^{-y}\right|x^{n}\right\rangle \nonumber \\
= & \left\langle \left.\partial_{t}\left(e^{-t}\lif_{k}\left(\log\left(1+\frac{t}{a}\right)\right)\left(1+\frac{t}{a}\right)^{-y}\right)\right|x^{n-1}\right\rangle \nonumber \\
= & \left\langle \left.\left(\partial_{t}e^{-t}\right)\lif_{k}\left(\log\left(1+\frac{t}{a}\right)\right)\left(1+\frac{t}{a}\right)^{-y}\right|x^{n-1}\right\rangle \nonumber \\
 & +\left\langle \left.e^{-t}\left(\partial_{t}\lif_{k}\left(\log\left(1+\frac{t}{a}\right)\right)\right)\left(1+\frac{t}{a}\right)^{-y}\right|x^{n-1}\right\rangle \nonumber \\
 & +\left\langle \left.e^{-t}\lif_{k}\left(\log\left(1+\frac{t}{a}\right)\right)\left(\partial_{t}\left(1+\frac{t}{a}\right)^{-y}\right)\right|x^{n-1}\right\rangle \nonumber \\
= & -\left\langle \left.e^{-t}\lif_{k}\left(\log\left(1+\frac{t}{a}\right)\right)\left(1+\frac{t}{a}\right)^{-y}\right|x^{n-1}\right\rangle \nonumber \\
 & +\left\langle \left.e^{-t}\left(\partial_{t}\lif_{k}\left(\log\left(1+\frac{t}{a}\right)\right)\right)\left(1+\frac{t}{a}\right)^{-y}\right|x^{n-1}\right\rangle \nonumber \\
 & -\frac{y}{a}\left\langle \left.e^{-t}\lif_{k}\left(\log\left(1+\frac{t}{a}\right)\right)\left(1+\frac{t}{a}\right)^{-y-1}\right|x^{n-1}\right\rangle \nonumber \\
= & -PC_{n-1}^{\left(k\right)}\left(y:a\right)-\frac{1}{a}yPC_{n-1}^{\left(k\right)}\left(y+1:a\right)\nonumber \\
 & +\left\langle \left.e^{-t}\left(\partial_{t}\lif_{k}\left(\log\left(1+\frac{t}{a}\right)\right)\right)\left(1+\frac{t}{a}\right)^{-y}\right|x^{n-1}\right\rangle .\nonumber
\end{align}

Now, we observe that
\begin{align}
 & \partial_{t}\lif_{k}\left(\log\left(1+\frac{t}{a}\right)\right)\label{eq:57}\\
= & \frac{1}{a\left(1+\frac{t}{a}\right)\log\left(1+\frac{t}{a}\right)}\left\{ \lif_{k-1}\left(\log\left(1+\frac{t}{a}\right)\right)-\lif_{k}\left(\log\left(1+\frac{t}{a}\right)\right)\right\} .\nonumber
\end{align}

From (\ref{eq:57}), we have
\begin{align}
 & \left\langle \left.e^{-t}\left(\partial_{t}\lif_{k}\left(\log\left(1+\frac{t}{a}\right)\right)\right)\left(1+\frac{t}{a}\right)^{-y}\right|x^{n-1}\right\rangle \label{eq:58}\\
= & \frac{1}{a}\left\langle \left.e^{-t}\frac{\lif_{k-1}\left(\log\left(1+\frac{t}{a}\right)\right)-\lif_{k}\left(\log\left(1+\frac{t}{a}\right)\right)}{\left(1+\frac{t}{a}\right)\log\left(1+\frac{t}{a}\right)}\left(1+\frac{t}{a}\right)^{-y}\right|x^{n-1}\right\rangle \nonumber \\
= & \frac{1}{n}\sum_{l=0}^{n-1}\dbinom{n}{l}\frac{\hat{C}_{l}}{a^{l}}\left\{ PC_{n-l}^{\left(k-1\right)}\left(y:a\right)-PC_{n-l}^{\left(k\right)}\left(y:a\right)\right\} .\nonumber
\end{align}

Therefore, by (\ref{eq:56}) and (\ref{eq:58}), we obtain the following
theorem.
\begin{thm}
\label{thm:6} For $n\ge0$, $k\in\mathbb{Z}$ and $a\ne0$, we have
\begin{eqnarray*}
PC_{n}^{\left(k\right)}\left(x:a\right) & = & -PC_{n-1}^{\left(k\right)}\left(x:a\right)-\frac{1}{a}xPC_{n-1}^{\left(k\right)}\left(x+1:a\right)\\
 & + & \frac{1}{n}\sum_{l=0}^{n-1}\dbinom{n}{l}\frac{\hat{C}_{l}}{a^{l}}\left\{ PC_{n-l}^{\left(k-1\right)}\left(x:a\right)-PC_{n-l}^{\left(k\right)}\left(x:a\right)\right\} .
\end{eqnarray*}
\end{thm}
\begin{rem*}
Note that
\begin{align}
 & \frac{1}{a}\left\langle e^{-t}\left(\frac{\lif_{k-1}\left(\log\left(1+\frac{t}{a}\right)\right)-\lif_{k}\left(\log\left(1+\frac{t}{a}\right)\right)}{\frac{t}{a}}\right)\right.\label{eq:59}\\
 & \qquad\qquad\qquad\qquad\qquad\times\left.\left(1+\frac{t}{a}\right)^{-y-1}\right|\left.\frac{\frac{t}{a}}{\log\left(1+\frac{t}{a}\right)}x^{n-1}\right\rangle \nonumber \\
= & \sum_{l=0}^{n-1}\frac{C_{l}}{a^{l}}\dbinom{n-1}{l}\nonumber \\
 & \times\left\langle \left.e^{-t}\left(\frac{\lif_{k-1}\left(\log\left(1+\frac{t}{a}\right)\right)-\lif_{k}\left(\log\left(1+\frac{t}{a}\right)\right)}{t}\right)\left(1+\frac{t}{a}\right)^{-y-1}\right|t\frac{x^{n-l}}{n-l}\right\rangle \nonumber \\
= & \sum_{l=0}^{n-1}\dbinom{n-1}{l}\frac{1}{n-l}\frac{C_{l}}{a^{l}}\left\langle e^{-t}\left(\lif_{k-1}\left(\log\left(1+\frac{t}{a}\right)\right)\right.\right.\nonumber \\
 & \qquad\qquad\qquad\qquad\quad\left.\left.\left.-\lif_{k}\left(\log\left(1+\frac{t}{a}\right)\right)\right)\left(1+\frac{t}{a}\right)^{-y-1}\right|x^{n-l}\right\rangle \nonumber \\
= & \frac{1}{n}\sum_{l=0}^{n-1}\dbinom{n}{l}\frac{C_{l}}{a^{l}}\nonumber \\
 & \times\left\langle \left.e^{-t}\left(\lif_{k-1}\left(\log\left(1+\frac{t}{a}\right)\right)-\lif_{k}\left(\log\left(1+\frac{t}{a}\right)\right)\right)\left(1+\frac{t}{a}\right)^{-y-1}\right|x^{n-l}\right\rangle \nonumber \\
= & \frac{1}{n}\sum_{l=0}^{n-1}\dbinom{n}{l}\frac{C_{l}}{a^{l}}\left\{ PC_{n-l}^{\left(k-1\right)}\left(y+1:a\right)-PC_{n-l}^{\left(k\right)}\left(y+1:a\right)\right\} .\nonumber
\end{align}

\end{rem*}
$\,$

By (\ref{eq:56}) and (\ref{eq:59}), we also get the following equation
:
\begin{align}
 & PC_{n}^{\left(k\right)}\left(x:a\right)\label{eq:60}\\
= & -PC_{n-1}^{\left(k\right)}\left(x:a\right)-\frac{1}{a}xPC_{n-1}^{\left(k\right)}\left(x+1:a\right)\nonumber \\
 & +\frac{1}{n}\sum_{l=0}^{n-1}\dbinom{n}{l}\frac{C_{l}}{a^{l}}\left\{ PC_{n-l}^{\left(k-1\right)}\left(x+1:a\right)-PC_{n-l}^{\left(k\right)}\left(x+1:a\right)\right\} .\nonumber
\end{align}

By the same method as Theorem \ref{thm:6}, we see that
\begin{align}
 & P\hat{C}_{n}^{\left(k\right)}\left(x:a\right)\label{eq:61}\\
= & -P\hat{C}_{n-1}^{\left(k\right)}\left(x:a\right)+\frac{1}{a}xP\hat{C}_{n-1}^{\left(k\right)}\left(x-1:a\right)\nonumber \\
 & +\frac{1}{n}\sum_{l=0}^{n-1}\frac{\hat{C}_{l}}{a^{l}}\dbinom{n}{l}\left\{ P\hat{C}_{n-l}^{\left(k-1\right)}\left(x:a\right)-P\hat{C}_{n-l}^{\left(k\right)}\left(x:a\right)\right\} ,\nonumber
\end{align}
and
\begin{align}
 & P\hat{C}_{n}^{\left(k\right)}\left(x:a\right)\label{eq:62}\\
= & -P\hat{C}_{n-1}^{\left(k\right)}\left(x:a\right)+\frac{1}{a}xP\hat{C}_{n-1}^{\left(k\right)}\left(x-1:a\right)\nonumber \\
 & +\frac{1}{n}\sum_{l=0}^{n-1}\dbinom{n}{l}\frac{C_{l}}{a^{l}}\left\{ P\hat{C}_{n-1}^{\left(k-1\right)}\left(x-1:a\right)-P\hat{C}_{n-l}^{\left(k\right)}\left(x-1:a\right)\right\} .\nonumber
\end{align}

Here, we compute
\[
\left\langle \left.e^{-t}\lif_{k}\left(-\log\left(1+\frac{t}{a}\right)\right)\left(\log\left(1+\frac{t}{a}\right)\right)^{m}\right|x^{n}\right\rangle
\]
in two different ways.

On the one hand,
\begin{align}
 & \left\langle \left.e^{-t}\lif_{k}\left(-\log\left(1+\frac{t}{a}\right)\right)\left(\log\left(1+\frac{t}{a}\right)^{m}\right)\right|x^{n}\right\rangle \label{eq:63}\\
= & \sum_{l=0}^{n-m}\frac{m!}{a^{l+m}}\dbinom{n}{l+m}S_{1}\left(l+m,m\right)\left\langle \left.e^{-t}\lif_{k}\left(-\log\left(1+\frac{t}{a}\right)\right)\right|x^{n-l-m}\right\rangle \nonumber \\
= & \sum_{l=0}^{n-m}\frac{m!}{a^{l+m}}\dbinom{n}{l+m}S_{1}\left(l+m,m\right)P\hat{C}_{n-l-m}^{\left(k\right)}\left(0:a\right)\nonumber \\
= & \sum_{l=0}^{n-m}\frac{m!}{a^{n-l}}\dbinom{n}{l}S_{1}\left(n-l,m\right)P\hat{C}_{l}^{\left(k\right)}\left(0:a\right).\nonumber
\end{align}

On the other hand,

\begin{align}
 & \left\langle \left.e^{-t}\lif_{k}\left(-\log\left(1+\frac{t}{a}\right)\right)\left(\log\left(1+\frac{t}{a}\right)\right)^{m}\right|x^{n}\right\rangle \label{eq:64}\\
= & \left\langle \left.e^{-t}\lif_{k}\left(-\log\left(1+\frac{t}{a}\right)\right)\left(\log\left(1+\frac{t}{a}\right)\right)^{m}\right|x\cdot x^{n-1}\right\rangle \nonumber \\
= & \left\langle \left.\partial_{t}\left\{ e^{-t}\lif_{k}\left(-\log\left(1+\frac{t}{a}\right)\right)\left(\log\left(1+\frac{t}{a}\right)\right)^{m}\right\} \right|x^{n-1}\right\rangle \nonumber \\
= & \left\langle \left.\left(\partial_{t}e^{-t}\right)\lif_{k}\left(-\log\left(1+\frac{t}{a}\right)\right)\left(\log\left(1+\frac{t}{a}\right)\right)^{m}\right|x^{n-1}\right\rangle \nonumber \\
 & +\left\langle \left.e^{-t}\left(\partial_{t}\lif_{k}\left(-\log\left(1+\frac{t}{a}\right)\right)\right)\left(\log\left(1+\frac{t}{a}\right)\right)^{m}\right|x^{n-1}\right\rangle \nonumber \\
 & +\left\langle \left.e^{-t}\lif_{k}\left(-\log\left(1+\frac{t}{a}\right)\right)\left(\partial_{t}\left(\log\left(1+\frac{t}{a}\right)\right)^{m}\right)\right|x^{n-1}\right\rangle \nonumber \\
= & -\left\langle \left.e^{-t}\lif_{k}\left(-\log\left(1+\frac{t}{a}\right)\right)\left(\log\left(1+\frac{t}{a}\right)\right)^{m}\right|x^{n-1}\right\rangle \nonumber \\
 & +\frac{m-1}{a}\left\langle \left.e^{-t}\lif_{k}\left(-\log\left(1+\frac{t}{a}\right)\right)\left(1+\frac{t}{a}\right)^{-1}\right|\left(\log\left(1+\frac{t}{a}\right)\right)^{m-1}x^{n-1}\right\rangle \nonumber \\
 & +\frac{1}{a}\left\langle \left.e^{-t}\lif_{k-1}\left(-\log\left(1+\frac{t}{a}\right)\right)\left(1+\frac{t}{a}\right)^{-1}\right|\left(\log\left(1+\frac{t}{a}\right)\right)^{m-1}x^{n-1}\right\rangle .\nonumber
\end{align}

It is easy to show that
\begin{align}
 & \left\langle \left.e^{-t}\lif_{k}\left(-\log\left(1+\frac{t}{a}\right)\right)\left(1+\frac{t}{a}\right)^{-1}\right|\log\left(1+\frac{t}{a}\right)^{m-1}x^{n-1}\right\rangle \label{eq:65}\\
= & \sum_{l=0}^{n-m}\frac{\left(m-1\right)!}{a^{l+m-1}}\dbinom{n-1}{l+m-1}S_{1}\left(l+m-1,m-1\right)P\hat{C}_{n-l-m}^{\left(k\right)}\left(-1:a\right)\nonumber \\
= & \sum_{l=0}^{n-m}\frac{\left(m-1\right)!}{a^{n-l-1}}\dbinom{n-1}{l}S_{1}\left(n-1-l,m-1\right)P\hat{C}_{l}^{\left(k\right)}\left(-1:a\right).\nonumber
\end{align}

Thus, by (\ref{eq:64}) and (\ref{eq:65}), we get
\begin{align}
 & \left\langle \left.e^{-t}\lif_{k}\left(-\log\left(1+\frac{t}{a}\right)\right)\left(\log\left(1+\frac{t}{a}\right)\right)^{m}\right|x^{n}\right\rangle \label{eq:66}\\
= & -\sum_{l=0}^{n-m-1}\frac{m!}{a^{n-l-1}}\dbinom{n-1}{l}S_{1}\left(n-1-l,m\right)P\hat{C}_{l}^{\left(k\right)}\left(0:a\right)\nonumber \\
 & +\left(\frac{m-1}{m}\right)\sum_{l=0}^{n-m}\frac{m!}{a^{n-l}}\dbinom{n-1}{l}S_{1}\left(n-l-1,m-1\right)P\hat{C}_{l}^{\left(k\right)}\left(-1:a\right)\nonumber \\
 & +\frac{1}{m}\sum_{l=0}^{n-m}\frac{m!}{a^{n-l}}\dbinom{n-1}{l}S_{1}\left(n-l-1,m-1\right)P\hat{C}_{l}^{\left(k-1\right)}\left(-1:a\right).\nonumber
\end{align}

Therefore, by (\ref{eq:63}) and (\ref{eq:66}), we obtain the following
theorem.
\begin{thm}
For $n,\, m\ge0$ with $n-m\ge0,$ $k\in\mathbb{Z}$ and $a\ne0$,
we have
\begin{align*}
 & \sum_{l=0}^{n-m}\frac{\dbinom{n}{l}}{a^{n-l}}S_{1}\left(n-l,m\right)P\hat{C}_{l}^{\left(k\right)}\left(0:a\right)\\
 & +\sum_{l=0}^{n-1-m}\frac{\dbinom{n-1}{l}}{a^{n-1-l}}S_{1}\left(n-1-l,m\right)P\hat{C}_{l}^{\left(k\right)}\left(0:a\right)\\
= & \left(1-\frac{1}{m}\right)\sum_{l=0}^{n-m}\frac{\dbinom{n-1}{l}}{a^{n-l}}S_{1}\left(n-1-l,m-1\right)P\hat{C}_{l}^{\left(k\right)}\left(-1:a\right)\\
 & +\frac{1}{m}\sum_{l=0}^{n-m}\frac{\dbinom{n-1}{l}}{a^{n-l}}S_{1}\left(n-1-l,m-1\right)P\hat{C}_{l}^{\left(k\right)}\left(-1:a\right).
\end{align*}
\end{thm}
\begin{rem*}
From the computation of
\begin{align*}
\left\langle \left.e^{-t}\lif_{k}\left(\log\left(1+\frac{t}{a}\right)\right)\left(\log\left(1+\frac{t}{a}\right)\right)^{m}\right|x^{n}\right\rangle
\end{align*}
 in two different ways, we can also derive the following equation
:

\begin{align}
 & \sum_{l=0}^{n-m}\frac{\dbinom{n}{l}}{a^{n-l}}S_{1}\left(n-l,m\right)PC_{l}^{\left(k\right)}\left(0:a\right)\label{eq:67}\\
 & +\sum_{l=0}^{n-1-m}\frac{\dbinom{n-1}{l}}{a^{n-l-1}}S_{1}\left(n-l-1,m\right)PC_{l}^{\left(k\right)}\left(0:a\right)\nonumber \\
= & \left(1-\frac{1}{m}\right)\sum_{l=0}^{n-m}\frac{\dbinom{n-1}{l}}{a^{n-l}}S_{1}\left(n-1-l,m-1\right)PC_{l}^{\left(k\right)}\left(1:a\right)\nonumber \\
 & +\frac{1}{m}\sum_{l=0}^{n-m}\frac{\dbinom{n-1}{l}}{a^{n-l}}S_{1}\left(n-1-l,m-1\right)PC_{l}^{\left(k-1\right)}\left(1:a\right).\nonumber
\end{align}

\end{rem*}
$\,$

By (\ref{eq:21}), (\ref{eq:26}) and (\ref{eq:27}), we easily see
that
\begin{equation}
\frac{d}{dx}PC_{n}^{\left(k\right)}\left(x:a\right)=\left(-1\right)^{n}n!\sum_{l=0}^{n-1}\frac{\left(-1\right)^{l}}{\left(n-l\right)l!a^{n-l}}PC_{l}^{\left(k\right)}\left(x:a\right),\label{eq:68}
\end{equation}

and
\begin{equation}
\frac{d}{dx}P\hat{C}_{n}^{\left(k\right)}\left(x:a\right)=\left(-1\right)^{n}n!\sum_{l=0}^{n-1}\frac{\left(-1\right)^{l-1}}{\left(n-l\right)l!a^{n-l}}P\hat{C}_{l}^{\left(k\right)}\left(x:a\right).\label{eq:69}
\end{equation}

For
\[
PC_{n}^{\left(k\right)}\left(x:a\right)\sim\left(e^{a\left(e^{-t}-1\right)}\frac{1}{\lif_{k}\left(-t\right)},a^{\left(e^{-t}-1\right)}\right),\;\left(a\ne0\right)
\]
and
\[
B_{n}^{\left(s\right)}\left(x\right)\sim\left(\left(\frac{e^{t}-1}{t}\right)^{s},t\right),\:\left(s\in\mathbb{Z}_{\ge0}\right),
\]
let us assume that
\begin{equation}
PC_{n}^{\left(k\right)}\left(x:a\right)=\sum_{m=0}^{n}C_{n,m}B_{m}^{\left(s\right)}\left(x\right).\label{eq:70}
\end{equation}

From (\ref{eq:25}), we note that
\begin{align}
 & C_{n,m}\label{eq:71}\\
= & \frac{\left(-1\right)^{m}}{m!}\nonumber \\
 & \times\left\langle \left.e^{-t}\lif_{k}\left(\log\left(1+\frac{t}{a}\right)\right)\left(1+\frac{t}{a}\right)^{-s}\left(\frac{\frac{t}{a}}{\log\left(1+\frac{t}{a}\right)}\right)^{s}\right|\left(\log\left(1+\frac{t}{a}\right)\right)^{m}x^{n}\right\rangle \nonumber
\end{align}

Now, we observe that
\begin{equation}
\left(\log\left(1+\frac{t}{a}\right)\right)^{m}x^{n}=\sum_{l=0}^{n-m}\frac{m!}{a^{n-l}}\dbinom{n}{l}S_{1}\left(n-l,m\right)x^{l}.\label{eq:72}
\end{equation}

By (\ref{eq:71}) and (\ref{eq:72}), we get
\begin{align}
 & C_{n,m}\label{eq:73}\\
= & \frac{\left(-1\right)^{m}}{m!}\sum_{l=0}^{n-m}\frac{m!}{a^{n-l}}\dbinom{n}{l}S_{1}\left(n-l,m\right)\nonumber \\
 & \qquad\times\left\langle \left.e^{-t}\lif_{k}\left(\log\left(1+\frac{t}{a}\right)\right)\left(1+\frac{t}{a}\right)^{-s}\right|\left(\frac{\frac{t}{a}}{\log\left(1+\frac{t}{a}\right)}\right)^{s}x^{l}\right\rangle \nonumber \\
= & \left(-1\right)^{m}\sum_{l=0}^{n-m}\frac{\dbinom{n}{l}}{a^{n-l}}S_{1}\left(n-l,m\right)\sum_{i=0}^{l}\frac{\dbinom{l}{i}\mathbb{C}_{i}^{\left(s\right)}}{a^{i}}\nonumber \\
 & \qquad\times\left\langle \left.e^{-t}\lif_{k}\left(\log\left(1+\frac{t}{a}\right)\right)\left(1+\frac{t}{a}\right)^{-s}\right|x^{l-i}\right\rangle \nonumber \\
= & \left(-1\right)^{m}\sum_{l=0}^{n-m}\frac{\dbinom{n}{l}}{a^{n-l}}S_{1}\left(n-l,m\right)\sum_{i=0}^{l}\frac{\dbinom{l}{i}\mathbb{C}_{i}^{\left(s\right)}}{a^{i}}PC_{l-i}^{\left(k\right)}\left(s:a\right)\nonumber \\
= & \left(-1\right)^{m}\sum_{l=0}^{n-m}\sum_{i=0}^{l}\frac{\dbinom{n}{l}\dbinom{l}{i}}{a^{n-l+i}}S_{1}\left(n-l,m\right)\mathbb{C}_{i}^{\left(s\right)}PC_{l-i}^{\left(k\right)}\left(s:a\right).\nonumber
\end{align}

Therefore, by (\ref{eq:70}) and (\ref{eq:73}), we obtain the following
theorem.
\begin{thm}
\label{thm:8} For $n\ge0$, $k\in\mathbb{Z}$ and $a\ne0$, we have
\begin{align*}
 & PC_{n}^{\left(k\right)}\left(x:a\right)\\
= & \sum_{m=0}^{n}\left\{ \left(-1\right)^{m}\sum_{l=0}^{n-m}\sum_{i=0}^{l}\frac{\dbinom{n}{l}\dbinom{l}{i}}{a^{n-l+i}}S_{1}\left(n-l,m\right)\mathbb{C}_{i}^{\left(s\right)}PC_{l-i}^{\left(k\right)}\left(s:a\right)\right\} B_{m}^{\left(s\right)}\left(x\right).
\end{align*}

\end{thm}
$\,$
\begin{rem*}
By the same method as Theorem \ref{thm:8}, we get
\begin{align}
 & P\hat{C}_{n}^{\left(k\right)}\left(x:a\right)\label{eq:74}\\
 & =\sum_{m=0}^{n}\left\{ \sum_{l=0}^{n-m}\sum_{i=0}^{l}\frac{\dbinom{n}{l}\dbinom{l}{i}}{a^{n-l+i}}S_{1}\left(n-l,m\right)\mathbb{\hat{C}}_{i}^{\left(s\right)}P\hat{C}_{l-i}^{\left(k\right)}\left(s:a\right)\right\} B_{m}^{\left(s\right)}\left(x\right).\nonumber
\end{align}

\end{rem*}
$\,$

For
\[
PC_{n}^{\left(k\right)}\left(x:a\right)\sim\left(e^{a\left(e^{-t}-1\right)}\frac{1}{\lif_{k}\left(-t\right)},a\left(e^{-t}-1\right)\right),\:\left(a\ne0\right),
\]

and
\[
H_{n}^{\left(s\right)}\left(x|\lambda\right)\sim\left(\left(\frac{e^{t}-\lambda}{1-\lambda}\right)^{s},t\right),\:\left(s\in\mathbb{Z}_{\ge0}\right),
\]
 let us assume that
\begin{equation}
PC_{n}^{\left(k\right)}\left(x:a\right)=\sum_{m=0}^{n}C_{n,m}H_{m}^{\left(s\right)}\left(x|\lambda\right).\label{eq:75}
\end{equation}

From (\ref{eq:25}), we have
\begin{align}
 & C_{n,m}\label{eq:76}\\
= & \frac{\left(-1\right)^{m}}{m!\left(1-\lambda\right)^{s}}\nonumber \\
 & \times\left\langle \left.e^{-t}\lif_{k}\left(\log\left(1+\frac{t}{a}\right)\right)\left(1+\frac{t}{a}\right)^{-s}\left(1-\lambda-\frac{\lambda t}{a}\right)^{-s}\right|\left(\log\left(1+\frac{t}{a}\right)\right)^{m}x^{n}\right\rangle \nonumber \\
= & \frac{\left(-1\right)^{m}}{\left(1-\lambda\right)^{s}}\sum_{l=0}^{n-m}\sum_{i=0}^{s}\frac{\dbinom{n}{l}\dbinom{s}{i}\left(l\right)_{i}}{a^{n-l+i}}\left(1-\lambda\right)^{s-i}\left(-\lambda\right)^{i}S_{1}\left(n-l,m\right)PC_{l-i}^{\left(k\right)}\left(s:a\right).\nonumber
\end{align}

Therefore, by (\ref{eq:75}) and (\ref{eq:76}), we obtain the following
theorem.
\begin{thm}
\label{thm:9} For $n\ge0$, $k\in\mathbb{Z}$ and $a\ne0$, we have
\begin{align*}
 & PC_{n}^{\left(k\right)}\left(x:a\right)\\
= & \frac{1}{\left(1-\lambda\right)^{s}}\sum_{m=0}^{n}\left\{ \left(-1\right)^{m}\sum_{l=0}^{n-m}\sum_{i=0}^{s}\frac{\dbinom{n}{l}\dbinom{l}{i}\left(l\right)_{i}}{a^{n-l+i}}\right.\\
 & \qquad\qquad\quad\left.\times\left(1-\lambda\right)^{s-i}\left(-\lambda\right)^{i}S_{1}\left(n-l,\, m\right)PC_{l-i}^{\left(k\right)}\left(s:a\right)\right\} H_{m}^{\left(s\right)}\left(x\mid\lambda\right).
\end{align*}

\end{thm}
$\,$

By the same method as Theorem \ref{thm:9}, we get {\small{
\begin{align}
 & P\hat{C}_{n}^{\left(k\right)}\left(x:a\right)\label{eq:77}\\
= & \frac{1}{\left(1-\lambda\right)^{s}}\sum_{m=0}^{n}\left\{ \sum_{l=0}^{n-m}\sum_{i=0}^{s}\frac{\dbinom{n}{l}\dbinom{s}{i}}{a^{n-l}}\left(-\lambda\right)^{s-i}S_{1}\left(n-l,m\right)P\hat{C}_{l}^{\left(k\right)}\left(i:a\right)\right\} H_{m}^{\left(s\right)}\left(x\mid\lambda\right).\nonumber
\end{align}
}}{\small \par}

For
\begin{eqnarray*}
PC_{n}^{\left(k\right)}\left(x:a\right) & \sim & \left(e^{a\left(e^{-t}-1\right)}\frac{1}{\lif_{k}\left(-t\right)},a\left(e^{-t}-1\right)\right),\\
x^{\left(n\right)} & = & x\left(x+1\right)\cdots\left(x+n-1\right)\sim\left(1,1-e^{-t}\right),
\end{eqnarray*}
let us assume that
\begin{equation}
PC_{n}^{\left(k\right)}\left(x:a\right)=\sum_{m=0}^{n}C_{n,m}x^{\left(m\right)}.\label{eq:78}
\end{equation}

From (\ref{eq:25}), we have
\begin{eqnarray}
C_{n,m} & = & \frac{1}{m!\left(-a\right)^{m}}\left\langle \left.e^{-t}\lif_{k}\left(\log\left(1+\frac{t}{a}\right)\right)\right|t^{m}x^{n}\right\rangle \label{eq:79}\\
 & = & \frac{1}{\left(-a\right)^{m}}\dbinom{n}{m}\left\langle \left.e^{-t}\lif_{k}\left(\log\left(1+\frac{t}{a}\right)\right)\right|x^{n-m}\right\rangle \nonumber \\
 & = & \frac{1}{\left(-a\right)^{m}}\dbinom{n}{m}PC_{n-m}^{\left(k\right)}\left(0:a\right).\nonumber
\end{eqnarray}

Therefore, by (\ref{eq:78}) and (\ref{eq:79}), we obtain the following
theorem.
\begin{thm}
\label{thm:10} For $n\ge0$, $k\in\mathbb{Z}$ and $a\ne0$, we have
\[
PC_{n}^{\left(k\right)}\left(x:a\right)=\sum_{m=0}^{n}\frac{\dbinom{n}{m}}{\left(-a\right)^{m}}PC_{n-m}^{\left(k\right)}\left(0:a\right)x^{\left(m\right)},
\]
where $x^{\left(m\right)}=x\left(x+1\right)\cdots\left(x+m-1\right).$\end{thm}
\begin{rem*}
By the same method as Theorem \ref{thm:10}, we get
\[
P\hat{C}_{n}^{\left(k\right)}\left(x:a\right)=\sum_{m=0}^{n}\frac{\dbinom{n}{m}}{a^{m}}P\hat{C}_{n-m}^{\left(k\right)}\left(0:a\right)\left(x\right)_{m},
\]
where $\left(x\right)_{m}=x\left(x-1\right)\cdots\left(x-m+1\right).$\end{rem*}


\bigskip
ACKNOWLEDGEMENTS. This work was supported by the National Research Foundation of Korea(NRF) grant funded by the Korea government(MOE)\\
(No.2012R1A1A2003786 ).
\bigskip


$\,$

\noindent Department of Mathematics, Sogang University, Seoul 121-742,
Republic of Korea

\noindent e-mail : dskim@sogang.ac.kr

\noindent $\,$

\noindent Department of Mathematics, Kwangwoon University, Seoul 139-701,
Republic of Korea

\noindent e-mail : tkkim@kw.ac.kr
\end{document}